\documentclass[11pt]{amsart}

\newcommand{\eps}{\epsilon}
\newcommand{\pa}{\partial}

\newfont{\fnt}{cmr10 scaled 550}
\renewcommand{\eps}{\varepsilon} 

\newtheorem{theorem}{Theorem}[section]
\newtheorem*{thm}{Theorem}

\newtheorem{lemma}{Lemma}[section]
\newtheorem{cor}{Corollary}[section]
\newtheorem{prop}{Proposition}[section] 

\newtheorem{definition}{Definition}[section]

\newtheorem{example}{Example}
\theoremstyle{remark}
\newtheorem{rem}{Remark}[section] 
\numberwithin{equation}{section}

\font\strange=msbm10

\renewcommand{\epsilon}{\varepsilon}

\renewcommand{\Sigma}{\varSigma}
\newcommand{\R}{{{\mathchoice  {\hbox{$\textstyle{\text{\strange R}}$}}
{\hbox{$\textstyle{\text{\strange R}}$}}
{\hbox{$\scriptstyle  N\kern-0.3em  R$}}  
{\hbox{$\scriptscriptstyle  R\kern-0.2em  R$}}}}}

\newcommand{\Z}{{{\mathchoice  {\hbox{$\textstyle{\text{\strange Z}}$}}
{\hbox{$\textstyle{\text{\strange Z}}$}}
{\hbox{$\scriptstyle  Z\kern-0.3em  Z$}}
{\hbox{$\scriptscriptstyle  Z\kern-0.2em  Z$}}}}}

\newcommand{\N}{{{\mathchoice  {\hbox{$\textstyle{\text{\strange N}}$}}
{\hbox{$\textstyle{\text{\strange N}}$}}
{\hbox{$\scriptstyle  N\kern-0.3em  N$}}
{\hbox{$\scriptscriptstyle  N\kern-0.2em  N$}}}}}

\renewcommand{\phi}{\varphi}
\usepackage{amsmath,amsthm,amscd}

\begin{document}
\title[Spectrum of quantum layers]{Existence of Bound States
for Layers Built Over Hypersurfaces in $\R^{n+1}$}
\date{January 6, 2004}
\author{Christopher Lin}
 \author{Zhiqin Lu}

\address{Department of
Mathematics, University of California,
Irvine, Irvine, CA 92697}

\email[Christopher Lin]{clin@math.uci.edu}
\email[Zhiqin Lu]{zlu@math.uci.edu}

\thanks{The first author is partially supported by NSF grant DMS-0202508.
The
second author is partially supported by  NSF Career award DMS-0347033 and the
Alfred P. Sloan Research Fellowship.}

\maketitle 

\tableofcontents 
\pagestyle{myheadings}

\newcommand{\ka}{K\"ahler }
\newcommand{\ii}{\sqrt{-1}}

\section{Introduction} 
 
In their study of the spectrum of quantum layers~\cite{DEK},
Duclos, Exner, and Krej\v{c}i\v{r}{\'\i}k proved the
existence
of bound states for certain quantum layers
\footnote{Quantum layers were studied by  many 
authors. An incomplete list of the works are
\cite{BEGK-1,CEK-1,DEK-1,EHL-1,ESK-2,
ESK-1,EK-1,PS-1,EV-1,kr-1,KA-1,KK-1}.
}. Part of their
motivations to study the quantum layers is from mesoscopic
physics. From the mathematical point of view, a quantum layer
is a 
 noncompact noncomplete manifold. For such a 
manifold, the spectrum of the Laplacian (with
Dirichlet or Neumann boundary condition) is less
understood. Nevertheless, from ~\cite{DEK} and this paper, we shall
see that the spectrum of a quantum layer has very interesting 
properties not only
from  the point of view of physics but also from the point of
view of mathematics. 

Mathematicians are  interested in the spectrum of two kinds
of manifolds: compact manifolds (with or without boundary),
and noncompact complete  manifolds.
For these two kinds of manifolds,
one can prove~\cite{Gaffney2,Gaffney1} that the 
Laplacians  can be uniquely extended as 
self-adjoint operators from operators on
smooth functions with compact support. For a compact manifold, by
the Hodge theorem, we can prove that the spectrum
of the Laplacian is composed of only discrete spectrum. On the
other hand, the spectrum of Laplacian of a 
noncompact complete manifold
is rather complicated. In general it has  
both discrete and essential spectrum.

In general, it is rather difficult to prove
the existence of discrete spectrum for a noncompact
 manifold, because  the existence of an
$L^2$ eigenfunction is a highly nontrivial fact.
However, in the following special  case, the discrete  spectrum does
exist.

We  define the following two
quantities:

\begin{definition}
Let $M$ be a manifold whose Laplacian $\Delta$ can be
extended to a self-adjoint operator. Let
\begin{align}
&\sigma_0=\inf_{f\in C_0^\infty(M)}\frac{-\int_M f\Delta
f}{\int_M f^2},\label{1-1}\\ &\sigma_{\rm ess}=
\sup_{K}\inf_{f\in C_0^\infty(M\backslash K)}
\frac{-\int_M f\Delta
f}{\int_M f^2},\label{1-2}
\end{align}
where $K$ is running over all compact subsets of $M$.
\end{definition}

$\sigma_0$ is the lower bound of the spectrum and 
$\sigma_{\rm ess}$ is the lower bound of the essential 
spectrum. In general, we have $\sigma_0\leq\sigma_{\rm
ess}$. If $\sigma_0<\sigma_{\rm ess}$, then the 
set of discrete spectrum must
be nonempty.  In particular, since the spectrum of a self-adjoint operator 
is a closed subset of the real line, there
is an
$L^2$ smooth function $f$ of $M$ such that
\[
\Delta f=-\sigma_0 f.
\]
$(\sigma_0,f)$ is called the ground state of the Laplacian.

We don't expect $\sigma_0<\sigma_{\rm ess}$ to be true 
in general. It seems that more often we would get
the opposite result $\sigma_0=\sigma_{\rm ess}$.
For example,
 by a theorem of Li and
Wang~\cite[Theorem 1.4]{Li-Wang}, we know  that if the
volume growth of a complete manifold is sub-exponential and
if the volume is infinite, then
$\sigma_{\rm ess}=0$. Thus in that case,
$\sigma_0=\sigma_{\rm ess}=0$ and we don't know a general 
way
to prove the existence of  ground state.

Let $\Sigma$
be an oriented $n$-manifold and $\pi:\Sigma\rightarrow
\R^{n+1}$ be an immersion. Let $N$ be the unit
normal vector field of $\Sigma$, we can define
the following map
\[
p:\Sigma\times(-a,a)\rightarrow \R^{n+1},\quad
(x,u)\mapsto x+uN,
\]
where $a$ is a small positive  number such that the map
$p$ is an immersion.

The quantum layer $\Omega$ built over $\Sigma$,
as a differentiable manifold, is very simple:
$\Omega=\Sigma\times (-a,a)$. The Riemannian metric
on $\Omega$ is defined by $p^*(ds_E^2)$, where $ds_E^2$
is the Euclidean metric of $\R^{n+1}$. The number
$a$ is called the depth of the layer.

\begin{rem}
The setting above is a little  bit more general than
in the paper~\cite{DEK}, where the authors require
that both $\Sigma$ and $\Omega$ are embedded. In particular,
they assume that the quantum layers are not
self-intersecting. There are some advantages of
our treatment: first, all the theorems still
remain true in immersed cases, and second, it is
possible to estimate the range of the depth $a$
using the upper bound of the second fundamental
form in the case of immersion, while in the
embedded case, global conditions of
$(\Sigma,\pi)$ must be imposed in order to keep
the layers from self-intersecting.
\end{rem}

The aim
of this paper is to study the ground state of the noncompact
noncomplete Riemannian manifold $(\Omega,p^*(ds_E^2))$,
where we assume  the Dirichlet boundary condition 
for the Laplacian. Our work is clearly motivated by the
work of 
~\cite{DEK}.

The first main result of this paper is the 
existence  of the ground state of 
layer over convex surface in $\R^3$. 
We are motivated by the following result in~\cite{DEK}:

\begin{thm}[Duclos, Exner, and Krej\v{c}i\v{r}{\'\i}k]
Let $\Omega$ be a layer of depth $a$ over a surface of
revolution whose Gauss curvature is integrable.
Suppose
$\Omega$ is not self-intersecting, and suppose
$a||A||<C_0<1$, where $||A||$ is the norm
of the second fundamental form and $C_0$ is a constant.
If the surface is not totally geodesic, then $\sigma_0
<\pi^2/(4a^2)$.
\end{thm}

Overlapping with the above result, we proved the following

 \begin{theorem}\label{main3}
Let $\Sigma$ be a convex surface in $\R^3$ which can be
represented by the graph of a  convex function
$z=f(x,y)$. Suppose $0$ is the minimum point of
the function and suppose that at $0$, $f$ is
strictly convex. Furthermore suppose that the
second fundamental form goes to zero at infinity.
Let $C$ be the supremum of the second fundamental
form of $\Sigma$. Let $Ca<1$. Then the ground
state of the quantum layer $\Omega$  exists.
\end{theorem}

\begin{rem}
We let $\Sigma$ to be the surface defined by the function
\[
f(x,y)=x^2+y^2.
\]
A straightforward computation gives the mean curvature of $\Sigma$:
\[
H=4\cdot\frac{1-\frac{2(x^2+y^2)}{1+4(x^2+y^2)}}{\sqrt{1+4(x^2+y^2)}}.
\]
Thus $H$ and then the second fundamental form goes to zero at infinity. By the above theorem, the quantum layer built from the above surface has a ground state.
\end{rem}

The second main result of this paper is motivated by the following:

\begin{thm}[Duclos, Exner, and Krej\v{c}i\v{r}{\'\i}k]
Let $\Sigma$ be a $C^2$-smooth complete simply connected
noncompact surface with a pole embedded in $\R^3$.
Let
the layer $\Omega$ built over the surface be not
self-intersecting. If the surface is not a  plane but it
is asymptotically planar, and if
the Gauss curvature is integrable and the total Gauss
curvature is nonpositive, then the ground state exists.  
\end{thm} 
  
\smallskip

\smallskip

In a more recent paper~\cite{CEK-1}, Carron, Exner, and 
Krej\v{c}i\v{r}{\'\i}k
 observed that  the
assumptions of simply-connectedness and the existence of a pole
on $\Sigma$
can be removed. 
$\Sigma$ is allowed to have a rather complicated topology.

\smallskip
\smallskip

By a
theorem of  Huber~\cite{Huber}, $\Sigma$ is conformally
equivalent to a compact Riemann surface with finitely many
points removed. In particular, we have
\[
\int_\Sigma K\leq 2\pi e(\Sigma),
\]
where $e(\Sigma)$ is the Euler characteristic
number of $\Sigma$.
The deficit of the above inequality can be
represented by isoperimetric constants. Let
$E_1,\cdots,E_s$ be the ends of the surface
$\Sigma$. For each  $E_i$, we define
\begin{equation}\label{iso}
\lambda_i=\underset{r\rightarrow\infty}{\lim}\,\frac
{A_i(r)}{\pi r^2},
\end{equation}
where $A_i(r)$ is the area of the ball $B(r)\cap
E_i$. We have the following 

\begin{thm}[Hartman~\cite{Hartman}] Using the
above notations, we have
\[
\frac{1}{2\pi}\int_\Sigma K=e(\Sigma)-\sum\lambda_i.
\]
\end{thm} 

\qed

We have the following
\footnote{There is an overlap of this result with the one in
~\cite{CEK-1}. The proofs are similar but not identical. 
In particular, we use the result of Hartman instead of the 
Kohn-Vossen inequality.}

\begin{theorem}\label{main2}
Suppose that $\Sigma$ is a complete immersed surface of
$\R^3$ such that the second fundamental form $A\rightarrow
0$. Suppose that the Gauss curvature is integrable and suppose
that 
\begin{equation}\label{1-4}
e(\Sigma)-\sum\lambda_i\leq 0,
\end{equation}
where $\lambda_i$ is the isoperimetric constant
at each end defined in~\eqref{iso}.
Let $a$ be a positive number such that $a||A||<C_0<1$.
If $\Sigma$ is not totally geodesic, then the ground
state of the quantum layer $\Omega$ exists. In particular, if
$e(\Sigma)\leq 0$, then the ground state exists.
\end{theorem}

More generally,
one of the key observation of our paper
is that in order to  generalize the results in~\cite{DEK}
to high dimensions, we must assume the parabolicity of
the hypersurface $\Sigma$. The parabolicity
of complete manifold was introduced  by Li and
Tam~\cite{Li-Tam}  (see also the survey papers
~\cite{Li,Li-1}).
A surface
with a pole and $L^1$ Gauss curvature is
parabolic. Thus the following result is a high
dimensional generalization of the  above
result of Duclos, Exner, and
Krej\v{c}i\v{r}{\'\i}k:

\begin{theorem}[{\bf Main theorem}]\label{main1}
Let $n\geq 2$ be a natural  number. Suppose $\Sigma\subset
\R^{n+1}$ is a complete immersed parabolic  hypersurface such that 
the second fundamental form $A\rightarrow 0$ at infinity.
Moreover, we assume that
\begin{equation}\label{133}
\sum_{k=1}^{[ n/2]}\mu_{2k}{\rm Tr}(\mathcal R^k)\,\, 
\text{is integrable and }\quad
\int_\Sigma \sum_{k=1}^{[n/2
]}\mu_{2k}{\rm Tr}(\mathcal R^k)
      d\varSigma \leq  0,
\end{equation}
where 
$\mu_{2k}>0$ for $k\geq 1$ are coefficients defined in
Lemma~\ref{lem51}, 
$[ n/2 ] $ is the integer part of $n/2$,
and ${\rm Tr}(\mathcal R^k)$ is defined in~\eqref{trace-1}. Let $a$ be a positive real number such that $a||A||<C_0<1$
for a constant $C_0$.
If $\Sigma$ is
not totally geodesic, then the ground state of the quantum
layer $\Omega$ exists.
\end{theorem}

\begin{cor} \label{main4}
Let $\rho$ be the scalar curvature of $\Sigma$. 
If $n=3$, then the main conditions~\eqref{133} in Theorem
~\ref{main1}
become
\begin{enumerate}
\item $\rho$ is integrable;
\item $\int_\Sigma\rho \,d\Sigma\leq 0$,
\end{enumerate}
If $n=4$, and if the sectional curvature of $\Sigma$
is positive outside a compact set of $\Sigma$, then the
conditions~\eqref{133} become
\begin{enumerate}
\item $\rho$ is
integrable;
\item $\int_\Sigma\rho \,d\Sigma
+16(\frac{\pi^2}{6}-1){a^3}
e(\Sigma)\leq 0$,
\end{enumerate}
where $e(\Sigma)$ is the Euler characteristic number of $\Sigma$.
\end{cor}

\smallskip

\smallskip

The organization of the paper is as follows: in \S 2, we define the quantum layers and give their basic properties; in \S 3, we give the lower bound of the essential spectrum of a quantum layer; in \S 4, the parabolicity of a submanifold of $\R^{n+1}$ is introduced; in \S 5, Theorem~\ref{main2}, The main theorem (Theorem~\ref{main1}), and Corollary~\ref{main4} are proved; finally, in \S 6, Theorem ~\ref{main3} is proved.

\smallskip 

We end up this section by posing the following
question: 

\begin{it}
Let $\Sigma$ be a noncompact complete Riemannian manifold of dimension
$n$. Then what  do we have to assume on $\Sigma$ so that when
$\Sigma\rightarrow \R^{n+1}$ is an asymptotically flat but not
totally geodesic immersion, the layer $\Omega$ built
over $\Sigma$ has ground state?

In particular, if $n=2$, the works of ~\cite{DEK,CEK-1} suggest that
the quantum layer $\Omega$ should have ground state
when the Gauss curvature is integrable\footnote{This part of the
question was implied in~\cite{DEK}. By ~\cite{DEK} and this paper,
we just need to show that for layers bulit over simply-connected
surfaces with positive total Gauss curvature, the ground state
exists.}.
\end{it}

\smallskip

\smallskip
 
{\bf Acknowledgement.} The authors thank P.~Li
for  advising us  to study the parabolicity
of manifolds, which makes this paper possible. The first
author thanks A.~Klein for clarifying a lot of the
subtleties of self-adjointness.
Corollary~\ref{main2} is the outcome of a
discussion
with G.~Tian. Finally, we give special thanks
to P. Exner and 
 D. Krej\v{c}i\v{r}{\'\i}k
for many useful suggestions to make the paper in its present form.

\section{Geometry of quantum layers}\label{S2}
Let $n>1$ be an integer and let $\Sigma$ be an immersed
(oriented)
hypersurface of $\R^{n+1}$. Let $a>0$ be a real number.
Heuristically 
speaking, the quantum layer $\Omega$ is obtained by
fattening
$\varSigma$ by a thickness of $a$ in the directions of $N$
and $-N$, respectively, where $N$ is the unit 
normal vector field. As a differentiable manifold,
$\Omega$ is just $\Sigma\times (-a,a)$.
We impose the following assumptions on
$\Sigma$ and $\Omega$:

\begin{enumerate}
\item[{\bf A1).}] Let $A$ be the second fundamental
form of $\Sigma$. We regard $A$ as  a linear operators
on $T_x\Sigma$ for every $x\in\Sigma$. We assume that there
is a constant $C_0$ such that $a||A||(x)<C_0<1$.

\item[{\bf A2).}] $\lVert A\rVert(x)
\rightarrow 0$ as $d(x,x_0)\rightarrow \infty$,
where $x_0\in \Sigma$ is a fixed  point.

\end{enumerate}

\newtheorem{remark}{Remark}[section]

\begin{definition}\label{def21}
 Let
$x_1,\cdots,x_n$ be a local coordinate system of 
$\Sigma$. Then $(\frac{\pa}{\pa x_1},\cdots,\frac{\pa}{\pa
x_n},\frac{\pa}{\pa u})$ is a local frame of $\Omega$,
where $u\in(-a,a)$. Such a local coordinate system
of $\Omega$
is referred as a standard coordinate system of $\Omega$
in this paper.
\end{definition}

We consider the map
\[
p: \Sigma\times (-a,a)\rightarrow \R^{n+1},
\quad (x,u)\mapsto x+uN.
\]

 Let $y_1,\cdots,y_{n+1}$
be the Euclidean coordinates of $\R^{n+1}$. 
Let 
\[
ds_E^2=dy_1^2+\cdots+ dy_{n+1}^2
\]
be the Euclidean metric of $\R^{n+1}$.
Let
$G_{ij}$ ($i,j=1,\cdots,n+1$) be defined by
\[
\sum_{i,j=1}^n G_{ij} dx_idx_j+
\sum_{i=1}^nG_{i,n+1}dx_idu
+\sum_{i=1}^nG_{n+1,i}dudx_i
+G_{n+1,n+1}dudu=p^*(ds_E^2).
\]
If $p$ is nonsingular
at a point, then the matrix $G_{ij}$ is positive definite
at that point.
 In order to express $G_{ij}$ in term of the
geometry of $\Sigma$, we introduce the following notations:

Let
$(h_{ij})$ ($i,j=1,\cdots n$) be the 
matrix representation of the second fundamental form $A$
with respect to the local frame $(\frac{\pa}{\pa x_1},
\cdots,\frac{\pa}{\pa x_n})$. Let
$g_{ij}dx_idx_j=p^*(ds_E^2)$ be the Riemannian metric of
$\Sigma$.
Let $h^\sigma_j=g^{\sigma i}h_{ij}$. Then a straightforward
computation gives (cf.~\cite{DEK}):

\begin{equation} \label{chris}
G_{ij} = \left\{
\begin{array}{ll}
(\delta_{i}^{\sigma} -
uh_{i}^{\sigma})(\delta_{\sigma}^{\rho} -
uh_{\sigma}^{\rho})g_{\rho j}& 1\leq i,j\leq n\\
0 & i {\text \,\, or \,\,} j=n+1\\
1& i=j=n+1
\end{array}
\right..
\end{equation}
In particular, we have
\begin{equation}\label{2-1}
 \det(G_{ij}) = (\det
(1-uA))^{2}\det(g_{ij}),
\end{equation}
and
\begin{align}
\begin{split}
&
\det(I-uA) = \prod_{i=1}^{n}(1-u\lambda_{i})\\
& = 1 -
u\sum_{i=1}^{n}\lambda_{i} + u^{2}\sum_{i<j}\lambda_{i}\lambda_{j} -
u^{3}\sum_{i<j<l}\lambda_{i}\lambda_{j}\lambda_{l} + \cdots +
(-1)^{n}u^{n}\prod_{i=1}^{n}\lambda_{i},
\end{split}
\end{align}
where $\lambda_1,\cdots,\lambda_n$ are the eigenvalues, or the 
principal curvatures of the second fundamental form
$A$. In a more intrinsic way, let $c_k(A)$ be the
$k$-th elementary polynomial of $A$. Then we have
\begin{equation}\label{2-3}
\det(I-uA) = \sum_{i=0}^{n}(-1)^{k}u^{k}c_{k}(A),
\end{equation}
where we define $c_{0}(A)= 1$. 

The following lemma is elementary but important:

\begin{lemma}\label{lem2-1}
Using the above notations and under Assumption A1), we
have
\[
(1-|u|\cdot||A||)^n\leq   |\det(I-uA)| \leq
(1+|u|\cdot||A||)^n.
\]
\end{lemma}

The proof is elementary and is omitted.

{\qed}

\begin{cor}\label{cor21}
We adopt  the above notations and 
Assumption A1). Then the  map $p$ is an immersion. In that
case,
$p^*(ds_E^2)$ is a Riemannian metric  on $\Omega$. Let
$d\Omega$ be the  measure defined by the metric and let
$dud\Sigma$ be the product measure on $\Omega$. Then we have
\begin{equation}\label{26}
(1-|u|||A||)^ndud\Sigma\leq d\Omega\leq
(1+|u|||A||)^ndud\Sigma.
\end{equation}
\end{cor}

{\bf Proof.} By Assumption A1), ~\eqref{2-1}
and Lemma~\ref{lem2-1}, we know that
$\det (G_{ij})>0$. Thus $p$ is nonsingular. ~\eqref{26}
follows from Lemma ~\ref{lem2-1} directly.

\qed

\begin{definition}
We define
the quantum layer to be the Riemannian manifold $(\Omega,
p^*(ds_E^2))$, where $ds_E^2$ is the standard Euclidean
metric
of $\R^{n+1}$. The real numbers $a$ and  $d=2a$ are called
the depth and the width of the quantum layer, respectively.
\end{definition}

The Laplacian $\Delta=\Delta_\Omega$ can be written as
\begin{equation}
\Delta=\frac{1}{\sqrt{\det(G_{kl})}}
\sum_{i,j=1}^n\frac{\pa}{\pa x_i}\left(
G^{ij}\sqrt{\det (G_{kl})}\frac{\pa}{\pa x_j}
\right)+\frac{1}{\sqrt{\det(G_{kl})}}\frac{\pa}{\pa u}
\left(\sqrt{\det (G_{kl})}\frac{\pa}{\pa u}\right),
\end{equation}
where $(x_1,\cdots,x_n,u)$ is the local coordinates defined
in Definition~\ref{def21}.
We have
\begin{equation}
(\Delta F,G)=(F,\Delta G) \quad \forall F, G \in 
C_{0}^{\infty}(\Omega),
\end{equation}
where $(\cdot,\cdot)$ is the $L^2$ inner product
\begin{equation}
(F,G)=\int_\Omega FG d\Omega.
\end{equation}
The norm $||F||$ is defined as
$\sqrt{(F,F)}$. If $F,G$ are differentiable, we  define
\begin{equation}
(\nabla F,\nabla G)=\int_\Omega\left(\sum_{i,j=1}^{n}G^{ij}
\frac{\pa F}{\pa x_i}\frac{\pa G}{\pa x_j}+
\frac{\pa F}{\pa u}\frac{\pa G}{\pa u}\right) d\Omega.
\end{equation}
Also, we define $||\nabla F||=\sqrt{(\nabla F,\nabla F)}$.

 In the case of
compact manifold or noncompact complete manifold, the
self-adjointness of the Laplacians is  classical
~\cite{Gaffney2,Gaffney1}. A quantum layer  is a noncompact
noncomplete manifold. For such a manifold, we still have

\begin{prop}
$\Delta$ can be extended as a self-adjoint 
operator.
\end{prop}

{\bf Proof.} According to ~\cite{Reed-Simon-1}, we define
the Hilbert space $H_1$ to be the closure of the 
space $C^\infty_{0}(\Omega)$ under the
norm
\[
||F||_{H_1}=\sqrt{||F||^2+||\nabla F||^2}.\]
We define the sesquilinear  form
\[
Q(F,G)=(\nabla F, \nabla G),
\]
for functions $F,G\in H_1$.
 By~\cite[Theorem VIII.15]
{Reed-Simon-1}, $Q$ is the quadratic form of a unique
self-adjoint operator. Such an
operator is an extension of $\Delta$,
which we still denote as  $\Delta$.
 Furthermore,
by the relation of $\Delta$ with the quadratic form,
we can verify that $\sigma_0$ and $\sigma_{\rm ess}$
in~\eqref{1-1},~\eqref{1-2} are the infimum of the
spectrum and the essential spectrum of $\Delta$,
respectively.

\qed

\section{Lower bound of the essential spectrum.}
The boundaries of $\Omega$ are $\Sigma\times \{{\pm}a\}$,
which are smooth manifolds. 
It is not hard to see that ~\eqref{1-2} can be written as
\[
\sigma_{\rm ess}
= \lim_{i\to \infty} \inf \{\frac{\int_\Omega \mid 
\nabla f\mid^{2}}
{\int_\Omega f^{2}} \mid f\in C_{0}^{\infty}
(\Omega \setminus K_{i})\},
\] 
where $ \left\{x_{0}\right\} \subset K_{1} \subset K_{2}
\subset \cdots$  is any compact exhaustion of
$\Omega$.  For example,
 we
can take 
\[
K_{i} =  \{ x + uN \mid x\in \overline{B_{x_o}(i)}
\subset \Sigma, 
u\in\left[\frac{-a(i-1)}{i}, \frac{a(i-1)}{i}\right] \}.
\]

To perform our estimate we need to obtain
 a lower bound for the Rayleigh quotient $\frac{\int \mid
\nabla f\mid^{2}}{\int f^{2}}, \forall f\in
C_{0}^{\infty}(\Omega \setminus K_{i})$ for a large enough
$i\in \mathbb{N}$.

We use the standard coordinate system $(x_1,\cdots,x_n,u)$
of Definition~\ref{def21}. Let        
$ f_j = \frac{\partial f}{\partial x_j}$, ($i=1,\cdots n)$
and $f_{n+1}=\frac{\pa f}{\pa u}$. 
Then

\[
|\nabla f|^{2}                                  
 = |f_{n+1}|^{2} + \sum_{k,l \ne n+1}G^{kl}f_kf_l,
\]
where $G^{ij}$ is the inverse of $G_{ij}$.
In particular, we have
\begin{equation}\label{234}
|\nabla f|^{2} \geq \left|\frac{\partial f}{\partial
u}\right|^{2}.
\end{equation}
Since $f=0$ on $\pa\Omega$, the Poincar\'e inequality gives
\begin{equation}\label{3-4}
\int_{-a}^{a}\left|\frac{\partial f}{\partial u}
\right|^{2}\, du 
 \geq\kappa_1^{2}\int_{-a}^{a}f(u)^{2}\, du,
\end{equation}
where $\kappa_1=\pi/2a$.

\begin{theorem}\label{thm31}
 Under Assumptions A1) and A2), 
we have $\sigma_{\rm ess} \geq \kappa_{1}^{2}$. 
\end{theorem}

{\bf Proof.}
We first observe that for arbitrary $\eps>0$, there is an 
$i$ large enough such that $||A||<\eps$ on
$\Sigma\backslash
K_i$. By  Corollary ~\ref{cor21}, we know that
\begin{equation}\label{3-3}
(1-a\eps)^n dud\Sigma\leq d\Omega\leq
(1+a\eps)^ndud\Sigma.
\end{equation}
Thus we have
\begin{equation}\label{new}
\int_\Omega f^2d\Omega\leq(1+a\eps)^n
\int_\Sigma\int_{-a}^af^2 dud\Sigma.
\end{equation}
On the other hand, by~\eqref{234}, ~\eqref{3-4}
and~\eqref{3-3}, we have
\begin{equation}\label{new1}
\int_\Omega|\nabla f|^2d\Omega\geq
(1-a\eps)^n\kappa_1^2\int_\Sigma\int^a_{-a}f^2dud\Sigma.
\end{equation}
Comparing ~\eqref{new} and~\eqref{new1}, we have
\[
\sigma_{\rm ess}\geq\frac{(1-a\eps)^n}{(1+a\eps)^n}
\kappa_1^2.
\]
Since $\eps$ is arbitrary, 
we get the conclusion of the theorem.

\qed

\section{Parabolicity of complete Riemannian manifolds}
Before giving the formal definition, we study the
following example.

Suppose $n>1$ is an integer. Let $R>0$ be a big number. We
are interested in the set of functions
\[
F(R)=\{f\in C_0^n(\R^n)| f\equiv 1\,\, \text 
{for }\,\,
|x|<R, \text{$f$ is rotationally symmetric.} \}.
\]

We have the following

\begin{example}\label{exp1}
If $n>2$, then for any $C>0$ there exists an
$R_{0}$ such that  for any $R>R_{0}$ we have 
\[
\int_{\mathbb{R}^n}|\nabla f|^2 > C
\]
for all $f\in F(R)$. If $n=2$, then for any $\epsilon > 0$ there exists
$R_{0}>0$ such that for any $R>R_{0}$, we can find
an $f_{R} \in F(R)$ for which
\[
\int_{\mathbb{R}^2}|\nabla f|^2 < \epsilon.
\]
\end{example}

{\bf Proof.} If $n>2$, then
\[
\int_{R}^\infty
1/r^{n-1}dr=\frac{1}{n-2}\cdot\frac{1}{R^{n-2}}.
\]
Thus we have
\[
\int_{\mathbb R^n}|\nabla f|^2\geq
(n-2)cR^{n-2}\int_{ R}^\infty r^{n-1}\left|\frac{\pa
f}{\pa r}
\right|^2dr
\int_{R}^\infty1/r^{n-1}dr\geq (n-2)cR^{n-2}
\rightarrow
+\infty
\]
by Cauchy inequality, where $c$ is the volume of the 
unit $(n-1)$-sphere.
However, for $n=2$, we let $f_R=\sigma_R(|x|)\in F(R)$, where
 $\sigma_R(t)$ is defined as
\[
\sigma_R(t)=      
\left\{
\begin{array}{ll}
1& t\leq R\\
(1-\frac{\log R}{R})^{-1}(\frac{\log R}{\log t}-
\frac{\log R}{R}) & R<t\leq e^R\\
0& t\geq e^R
\end{array}
\right.  .  
\]                
A straightforward computation gives

\[
\int_0^\infty t|\sigma_R'(t)|^2 dt\leq \frac 43\,\frac{1}{\log
R} \quad \text{for}\, R>3,
\]
and thus 
\begin{equation}\label{4-1}
\int_{\mathbb R^2}|\nabla f_R|^2\rightarrow 0,\quad
R\rightarrow
\infty.
\end{equation}
This completes the proof. 

\qed

The phenomenon in  the above example 
can be explained by the result
of Cheng-Yau~\cite[Section 1]{CY2}.
In~\cite[Definition 0.3]{Li-Tam}, the authors defined the following

\begin{definition}
A complete manifold is said to be parabolic, if it does not
admit a positive Green's function. Otherwise it is said to be
nonparabolic.
\end{definition}

\begin{rem} 
According to this definition, $\R^n$ is a parabolic
manifold if and only if $n=2$. In particular,
~\eqref{4-1} follows from Proposition~\ref{prop41} below,
which is a result  given in ~\cite{Li-Tam}.
\end{rem}

\begin{prop}\label{prop41}
Let $\Sigma$ be a parabolic manifold.
Let $B(r)$ be the ball of radius $r$ in $\Sigma$ with respect to 
a reference point $x_0$. Let $R>r>1$. Consider the following
Dirichlet problem
\[
\left\{
\begin{array}{ll}
\Delta f=0,&{\text on\,\,} B(R)\backslash B(r)\\
f=0, & {\text on\,\,} \Sigma\backslash B(R)\\
f=1, & {\text on\,\,} B(r).
\end{array}
\right.
\]
Then we have
\[
\lim_{R\rightarrow \infty}\int_\Sigma|\nabla f|^2=0.
\]
\end{prop} 

\qed

\begin{rem} The functions $f$ serve as the 
high dimensional generalization of the MacDonald
functions in the paper~\cite[page 21]{DEK}.
 These functions
will play an important role in the next section. 
\end{rem}

The following  geometric criterion of parabolicity
was proved by Grigor'yan~\cite{Gro1,Gro2} 
and Varopoulos~\cite{varo} independently (cf.~\cite[Equation
(3.1)]{Li}):

\begin{theorem}\label{thm4111}
Let $V(t)$ be the volume of the geodesic ball $B(t)$.
If $\Sigma$ is nonparabolic, then 
\[
\int_1^\infty\frac{tdt}{V(t)}<\infty.
\]
In particular,
if $V(t)$ is  
at most of quadratic growth, then $\Sigma$ is parabolic.
\end{theorem}

\qed

\begin{cor}\label{cor41}
Let $\varSigma$
 be a smooth surface whose  Gauss curvature $K \in
L^{1}(\varSigma)$.  Then 
$\varSigma$ is  a parabolic manifold of
 dimension $2$.  
\end{cor}

{\bf Proof.}  We wish to
compare the volume  growth rate of the geodesic ball
$V(t)$ with $t$.  To do so, first we assume that $\Sigma$ has a pole
and we use  the polar coordinate system given by the exponential map
centered at
the pole  to  write 

 $$V(t) = \int_{0}^{t} \int_{0}^{2\pi} f(r,\theta)\, dr
d\theta,$$ 
where under  the polar  coordinates, the expression of the
metric becomes 
$ds_{\varSigma}^{2} = dr^{2} + f^{2}(r,\theta)d\theta^{2}$ on $\varSigma$.

It follows that
$$V'(t) = \int_{0}^{2\pi}
f(t,\theta)\, d\theta.$$ The Jacobi equation for the
exponential map gives 
$$f'' + Kf = 0 \hskip 0.1cm; \hskip 0.3cm f(0,\theta) = 0,
\hskip 0.1cm f'(0,\theta) = 1,$$
where the prime  denotes derivative in the radial direction. 
Thus we have
\[
V'''(t)=-\int_0^{2\pi}Kf(t,\theta) d\theta.
\]
Since $K$ is integrable, this implies that
\[
|V''(t)|\leq C
\]
for some constant $C$. Consequently,
\begin{equation}\label{poiq}
V(t)\leq Ct^2
\end{equation}
for $t$ large enough.

If the surface $\Sigma$ doesn't have
a pole, we  get the 
similar  estimate outside the cut locus with respect to a fixed
reference point. Since the measure of the cut locus is zero,
we get the same estimate~\eqref{poiq}. This is an observation of
Gromov.

Thus the volume of $\Sigma$ is at most of quadratic growth and it must
be parabolic by Theorem~\ref{thm4111}.

\qed

\section{The upper bound estimate of $\sigma_0$.} 
The idea to estimate $\sigma_0$,
the infimum of the spectrum of the Laplacian
from above is to
construct test functions which would provide the
strict upper bound
$\kappa_{1}^{2}$ (where $\kappa_1=\pi/2a$). 
We may construct test functions which are continuous
everywhere on
$\Omega$ and smooth everywhere on $\Omega$ except on a set
of measure 0.  
Such  functions must be in $H^{1,2}(\Omega)$,
which serve our purpose.

We define the quadratic form
\begin{equation}\label{5-1-1}
Q(\xi,\xi)=\int_\Omega|\nabla \xi|^2d\Omega
-\kappa_1^2\int_\Omega \xi^2d\Omega,
\end{equation}
for $\xi\in H^{1,2}(\Omega)$.
By the nature of the metric on $\Omega$, we define
\begin{equation}\label{5-2-1}
Q_1(\xi,\xi)=\int_\Omega|\nabla'\xi|^2d\Omega,
\end{equation}
where 
\begin{equation}\label{5-2-2}
|\nabla'\xi|^2=\sum_{i,j=1}^n
G^{ij}\frac{\pa\xi}{\pa
x_i}
\frac{\pa\xi}{\pa x_j},
\end{equation}
 and
\begin{equation}\label{5-3}
Q_2(\xi,\xi)=\int_\Omega\left(\frac{\pa \xi}{\pa u}\right)^2
d\Omega-\kappa_1^2\int_\Omega\xi^2d\Omega,
\end{equation}
where $(x_1,\cdots,x_n,u)$ are the standard
coordinates in Definition~\ref{def21}.
It is clear that 
\[
Q(\xi,\xi)=Q_1(\xi,\xi)+Q_2(\xi,\xi).
\]

The test functions we shall construct will essentially 
be the product of a vertical function (depending only on
$u$) and a horizontal one (depending only on
$x\in\varSigma$). Let $\phi=\psi\chi$ be a test function,
where $\psi\in C_0^\infty (\Sigma)$ and $\chi$ is a
smooth 
function of $u$ such that $\chi(\pm a)=0$.

  Note that
\[
\bigtriangledown (\chi\psi) = \chi\bigtriangledown \psi +
\psi\bigtriangledown \chi.
\]
By~\eqref{chris}, we have $<\nabla\psi,\nabla\chi>=0$.
Thus we have
\begin{equation} \label{5-1}
\int_{\Omega}|\bigtriangledown(\chi\psi)|^{2}=
\int_{\Omega}\chi^{2}|\bigtriangledown \psi |^{2} 
       + \int_{\Omega}\psi^{2}|\bigtriangledown \chi |^{2}.
\end{equation}   
We wish to prove, with the suitable choice of $\psi$ and
$\chi$, that
\begin{equation}\label{5-2}
Q(\phi,\phi) =
 \int_\Omega \phi_{u}^{2}    
 -  
\kappa_{1}^{2}\int_\Omega \phi^{2}
+\int_\Omega\chi^2|\nabla\psi|^2
<0,
\end{equation}

\noindent
where $\phi_{u}$ denotes $\frac{\partial \phi}{\partial u}$
and $\phi_i=\frac{\pa \phi}{\pa x_i}$
for $i=1,\cdots,n$. 

Like in the  paper~\cite{DEK}, we choose $\chi=\cos
\kappa_1u$.
We need the following elementary lemma:

\begin{lemma} \label{lem51}
Let $a>0$ be a positive number
and let $\kappa_1=\frac{\pi}{2a}$. Let $\chi(u) =
\cos{\kappa_1u}$,  let 
\begin{equation}
\mu_k=
\int_{-a}^{a} u^{k}(\chi_{u}^{2} -
\kappa_{1}^{2}\chi^{2})\, du 
 ,\, \forall k \geq 0.
\end{equation}
Then
\begin{equation}
\mu_k=\left\{
\begin{array}{ll}
0& \text{if $k$ is odd, or $k=0$};\\
\frac 12\frac{(k)!}{(2\kappa_1)^{k-1}}\sum_{l=1}^{k/2}
\frac{(-1)^{k/2-l}\pi^{2l-1}}{(2l-1)!},& \text{if $k\neq 0$
is even}.
\end{array}
\right.
\end{equation}
Furthermore, $\mu_k>0$ if $k\neq 0$ is even.
\end{lemma}

\qed

\begin{theorem}\label{thm51}
We assume that the hypersurface $\varSigma \subset
\mathbb{R}^{n+1}$ is parabolic satisfying Assumptions A1), A2). 
Moreover, we assume that
$\sum_{k=1}^{[ n/2]}\mu_{2k}c_{2k}(A)$  is integrable and  
\begin{equation}
\int_\Sigma \sum_{k=1}^{[n/2
]}\mu_{2k}c_{2k}(A)
      d\varSigma \leq  0,
\end{equation}
where $A$ is the second fundamental form of $\Sigma$,
$\mu_k$ for $k\geq 1$ is defined in
Lemma~\ref{lem51}, 
$[ n/2 ] $ is the integer part of $n/2$,
and $c_k(A)$ is the $k$-th elementary
symmetric polynomial of $A$. 
If $\Sigma$ is not totally geodesic, then
\[
\sigma_0 < \kappa_{1}^{2}
.\]
\end{theorem}

{\bf Proof.}
We first
consider the test functions of the form $\phi =
\psi\cdot\chi$.
We define $\psi$ as follows: Let $x_0$ be a fixed point of
$\Sigma$ and let $R>r>1$. Let $B(R)$ and $B(r)$ be two 
balls in $\Sigma$ of radius $R$ and $r$ 
centered at $x_0$
respectively. We define $\psi$ as
\begin{equation}
\left\{
\begin{array}{ll}
\Delta\psi=0 &\text{on $B(R)-B(r)$};\\
\psi|_{B(r)}\equiv 1;\\
\psi|_{\Sigma-B(R)}\equiv 0,
\end{array}
\right.
\end{equation}
and we define $\chi=\cos\kappa_1 u$.

By the definition of the functions $\chi$ and $\psi$, using
Lemma~\ref{lem2-1}, we know that there is a constant $C$
such that
 \begin{equation}\label{5-8}
       \int_{\Omega}\chi^{2}|\bigtriangledown \psi |^{2} 
d\Omega
 \leq C\int_\Sigma|\nabla_\Sigma\psi|^2d\Sigma,
\end{equation}
where $\nabla_\Sigma$ is the connection of
$\Sigma$.
We first assume that
\begin{equation}\label{5-9}
\int_\Sigma \sum_{k=1}^{[n/2
]}\mu_{2k}c_{2k}(A)
      d\varSigma =-\delta<0.
\end{equation}
By ~\eqref{2-3} and Lemma~\ref{lem51}, we have
\begin{equation}\label{5-10}
\int_\Omega \psi^{2}|\chi_u|^{2} - 
\kappa_{1}^{2}\int_\Omega \psi^{2}\chi^{2}
= \int_\Sigma \psi^{2}\sum_{k=1}^{[ n/2
]}\mu_{2k}c_{2k}(A)\, d\varSigma.
\end{equation}

By the maximum
 principle and the fact that $\psi|_{B(r)}
\equiv 1$,
we have

\begin{equation}
\int_\Sigma\psi^{2}\sum_{k\geq 1} \mu_{2k}c_{2k}(A) 
 \leq \int_{B(r)}\sum_{k\geq 1}\mu_{2k}c_{2k}(A) 
        +\int_{\varSigma \setminus B(r)}|
\sum_{k\geq 1}
\mu_{2k}c_{2k}(A)|.
\end{equation}

On the other side, since $\sum\mu_{2k}c_{2k}(A)$ is
integrable, if $r$ is large enough, by the above inequality,
we have
\begin{equation}\label{5-12}
\int_{\varSigma}\psi^2\sum
\mu_{2k}c_{2k}(A)d\Sigma< 
   -\frac{\delta}{2}.
\end{equation}
By Proposition~\ref{prop41} and~\eqref{5-8}, if $R$ large
enough, we have
\begin{equation}\label{5-13}
\int_\Omega\chi^2|\nabla\psi|^2\leq \frac{\delta}{4}.
\end{equation}
Combining~\eqref{5-10},~\eqref{5-12} and ~\eqref{5-13},
we proved ~\eqref{5-2} under the assumption~\eqref{5-9}.

Now we assume that

\begin{equation}\label{5-14}
\int_\Sigma \sum_{k=1}^{[n/2
]}\mu_{2k}c_{2k}(A)
      d\varSigma =0.
\end{equation}
In this case, the test functions $\phi=\psi\chi$ are not
good
enough to give the upper bound of $\sigma_0$. We shall use a
trick in~\cite{DEK} (see also~\cite{DP,RB}) to construct the
test functions. 

 We let
\[
\phi_\eps=\phi+\eps j\chi_1,
\]
where $\eps$ is a small number, $j$ is a smooth
function on
$\Sigma$ whose support is contained is $B(r-1)$,
and $\chi_1$ is a smooth function on $[-a,a]$ such that
$\chi_1(\pm a)=0$. As a general fact, we have
\begin{equation}\label{5-15}
Q(\phi_\eps,\phi_\eps)=Q(\phi,\phi)+2\eps Q(\phi,j\chi_1)
+\eps^2Q(j\chi_1,j\chi_1).
\end{equation}
 By~\eqref{5-2},~\eqref{5-8}, and~\eqref{5-10}, we
have
\begin{equation}\label{5-16}
Q(\phi_\eps,\phi_\eps)\leq C\int_\Sigma|\nabla\psi|^2
d\Sigma+
\int_\Sigma \psi^{2}\sum_{k=1}^{[ n/2
]}\mu_{2k}c_{2k}(A)\, d\varSigma+2\eps Q(\phi,j\chi_1)
+\eps^2Q(j\chi_1,j\chi_1).
\end{equation}
Since ${\rm supp}\, j\subset B(r-1)$, we have
\begin{equation}
Q(\phi,j\chi_1)=\int_\Omega j(\chi_u(\chi_1)_u-\kappa_1^2
\chi\chi_1)d\Omega=\int_\Sigma j
\int_{-a}^a(\chi_u(\chi_1)_u-\kappa_1^2
\chi\chi_1)\det(1-uA)dud\Sigma.
\end{equation}

Using integration by parts, we have
\begin{equation}
Q(\phi,j\chi_1)=-\int_\Sigma j\int_{-a}^a
\chi_u\frac{\pa}{\pa u}
\det (1-uA)\chi_1dud\Sigma.
\end{equation}

Now we are able to choose suitable $j$ and $\chi_1$ for our
purpose. By assumption, we know that $\Sigma$ is not
totally geodesic. Thus at least there is a point $x\in
\Sigma$ such that $\pa_u\det(1-uA)\not\equiv 0$. We assume that
$x\in B(r-1)$ without losing generality. We choose $\chi_1$
and $j$ such that the  integral $Q(\phi,j\chi_1)$  is not
zero. Note that the choice of $j$ is independent of $\phi$. We  choose 
$\eps$ (positive or negative) small enough so that 
\[
2\eps Q(\phi,j\chi_1)
+\eps^2Q(j\chi_1,j\chi_1)<0.
\]
Finally, since
\[
{\rm supp}\, j\subset B(r-1),
\]
 the above expression is independent to  
 $r$ and $R$. By the parabolicity of $\Sigma$, if $r,R\rightarrow \infty$, then
\[
\int_\Sigma|\nabla\psi|^2d\Sigma\rightarrow 0,
\]
and by the assumption~\eqref{5-14},
\[
\int_\Sigma \psi^{2}\sum_{k=1}^{[ n/2
]}\mu_{2k}c_{2k}(A)\, d\varSigma\rightarrow 0.
\]
Thus by~\eqref{5-16}, $Q(\phi_\eps,\phi_\eps)$ is negative for $r, R$ large.
This completes the proof of the theorem.

\qed

Let $\mathcal R=(R_{ijkl})$ be the curvature tensor of $\Sigma$. Define
\begin{equation}\label{trace-1}
{\rm tr}(\mathcal R^p)
=\sum_{i_s<j_s, k_s<l_s, 
s=1,\cdots,p} (-1)^{sgn(\sigma)}
R_{i_1j_1k_1l_1}\cdots R_{i_pj_pk_pl_p},
\end{equation}
where $\sigma$ is the permutation $(i_1,\cdots,j_p; k_1,\cdots,l_p)$.
Then from Gray~\cite[(4.15)]{gray}, we have 

\begin{prop}\label{lem63}
Using the above notations, we have
\[
{\rm Tr}(\mathcal R^p)=c_{2p}(A)
\]
\end{prop}

\qed

\begin{rem} If $n$ is even, then up to a constant,  ${\rm Tr}(\mathcal
R^{n/2})=c_n(A)$ is the Gauss-Bonnet-Chern density. 
\end{rem}

{\bf Proof of Theorem~\ref{main1}.}
By Theorem~\ref{thm31}, Theorem~\ref{thm51}, and Proposition~\ref{lem63}, we
have
\[
\sigma_0<\kappa_1^2\leq\sigma_{\rm ess}.
\]
Thus the ground state exists.

\qed 

{\bf Proof of Theorem~\ref{main2}.} By the Theorem of
Hartman, we know that ~\eqref{1-4} is equivalent
to 
\[
\int_\Sigma K\leq 0.
\]
Thus the result follows from Theorem~\ref{main1} for
$n=2$.

\qed

{\bf Proof of Corollary~\ref{main4}.} If $n=3$, then the conditions
~\eqref{133} are
\[
{\rm Tr}\,(\mathcal R^1) \text{ is integrable and}\quad
\int_\Sigma {\rm Tr}\,(\mathcal R^1)\leq 0.
\]
But $\rho=2\,{\rm Tr}(\mathcal R^1)$.

If $n=4$, a tedious  computation gives
\[
{\rm Tr}\, (\mathcal R^2)=\frac{1}{24}(\rho^2-4|{\rm
Ric}|^2+|{\mathcal  R}|^2), 
\]
where ${\rm Ric}$ is the Ricci curvature of $\Sigma$,
and $|{\rm Ric}|, |\mathcal R|$ are the norms of the Ricci tensor
and the curvature tensor, respectively. If the sectional curvature
is positive outside a compact set, then by~\cite[Theorem
9]{Greene-Wu},
${\rm Tr}\,(\mathcal R^2)$ is integrable and 
\[
\int_\Sigma {\rm Tr}\,(\mathcal R^2)\leq \frac{4\pi^2}{3}
e(\Sigma),
\]
where $e(\Sigma)$ is the Euler characteristic 
number of $\Sigma$.
The theorem follows from the above inequality, Proposition~\ref{lem51}, and
Theorem~\ref{main1}.

\qed

Before finishing this section, we give 
the following  example of  the manifold $\Sigma$ of
dimension $3$ satisfying the
conditions in Theorem~\ref{main1}. Thus the theorem is not an empty statement for high dimensions.

\begin{example}
There is a complete manifold $\Sigma$ of dimension
$3$ immersed in $\R^4$ such that
\begin{enumerate}
\item It is parabolic;
\item $A\rightarrow 0$, where $A$ is the
second fundamental form;
\item $\frac 12\int_\Sigma|\rho|=\int_\Sigma|c_2(A)|<+\infty$;
\item $\frac 12\int_\Sigma\rho=\int_\Sigma c_2(A)<0$.
\end{enumerate}
\end{example}

{\bf Proof.}
Let $\Sigma=S^1\times \R^2$.
We consider the immersion by
\[
\Sigma\rightarrow\R^4, \quad 
(\theta,t,\phi)\rightarrow 
(\sigma(t)\cos\theta,\sigma(t)\sin\theta,
t\cos\phi,t\sin\phi),
\]
where $\sigma(t)$ is a smooth positive function
defined in~\eqref{ex-n-1}.
Here we use $\theta$ as the local coordinate
of $S^1$ and $(x,y)\in \R^2$ with
$x=t\cos\phi, y=t\sin\phi$. 
The Riemannian metric of $\Sigma$ is
\[
ds^2=(1+\sigma'(t)^2)(dt)^2+\sigma^2(t)(d\theta)^2+t^2
(d\phi)^2.
\]
We claim that $\Sigma$ is parabolic. In order to 
prove this, we let $x_0=(1,0,0)\in \Sigma$. Let
$B(R)$ be the geodesic ball of radius $R$
centered at $x_0$.
Then $B(R)\subset \{x\in \varSigma | t<R\}$. To see this,
let $x\in B(R)$ such that ${\rm
dist}\,(x,x_0)=R'$, and let
$\eta=(\eta_1(s),\eta_2(s),\eta_3(s))$ be the
geodesic line of $\Sigma$ connecting $x_0$ and $x$,
where $s$ is the arc length. Then we
have
\[
R=R'\geq\int_0^{R'}(1+\sigma'(s)^2)^{\frac 12}
|\eta_1'(s)|ds\geq  t.
\]
From the above equation, we have
\[
{\rm vol}\,B(R)\leq 4\pi^2\int_0^{R}t\sigma(t)
(1+\sigma'(t)^2)^{\frac 12}dt\leq CR^2\log R
\]
for some constant $C$. Thus we have
\[
\int_0^\infty\frac{t}{{\rm vol}\,B(t)}dt
=+\infty,
\]
and  $\Sigma$ is parabolic by Theorem~\ref{thm4111}.

The normal vector of $\Sigma$ in $\R^4$ is
\[
N=\frac{1}{\sqrt{1+\sigma'(t)^2}}
(\cos\theta,\sin\theta,-\sigma'\cos\phi,-\sigma'\sin\phi
).
\]
The principal curvatures are
\[
\frac{\sigma''}{({1+\sigma'(t)^2})^{\frac
32}},\quad -\frac{1}{\sigma\sqrt{1+\sigma'(t)^2}},
\quad\frac{\sigma'}{t\sqrt{1+\sigma'(t)^2}}.
\]
By the definition of $\sigma(t)$, all principal
curvatures go to zero as $t\rightarrow \infty$. 
Thus $A\rightarrow 0$ at infinity. On the other 
hand
\begin{equation}\label{ex-n-2}
\int_\Sigma c_2(A)=4\pi^2
\int_0^\infty\left(\frac{\sigma\sigma'\sigma''}
{(1+\sigma'(t)^2)^{\frac 32}}
-\frac{t\sigma''}{(1+\sigma'(t)^2)^\frac 32}
-\frac{\sigma'}{\sqrt{1+\sigma'(t)^2}}
\right) dt
\end{equation}

 We let the function
$\sigma(t)$ be a smooth increasing function such
that
\begin{equation}\label{ex-n-1}
\left\{
\begin{array}{ll}
\sigma(t)=\log t& t>3+\eps\\
\sigma(t)=\log 3 & t<3
\end{array}
\right.,
\end{equation}
for $\eps$ small.
The last two terms of ~\eqref{ex-n-2} can be
calculated easily:
\begin{equation}\label{ex-n-3}
\int_0^\infty\left(
-\frac{t\sigma''}{(1+\sigma'(t)^2)^\frac 32}
-\frac{\sigma'}{\sqrt{1+\sigma'(t)^2}}
\right) dt=\left.-\frac{t\sigma'}
{\sqrt{1+\sigma'(t)^2}}\right|_0^\infty=-1.
\end{equation}
Let $R$ be a big number. We
have
\[
\int_0^R\frac{\sigma\sigma'\sigma''}
{(1+\sigma'(t)^2)^{\frac 32}}
dt
=-\frac{\log R}{(1+\frac{1}{R^2})^\frac 12}
+\log 3 +\int_{3+\eps}^R\frac{\sigma'(t)}
{(1+\sigma'(t)^2)^{\frac 12}} dt
+\int^{3+\eps}_3\frac{\sigma'(t)}
{(1+\sigma'(t)^2)^{\frac 12}} dt.
\]
The last term can be estimated by
\[
\int^{3+\eps}_3\frac{\sigma'(t)}
{(1+\sigma'(t)^2)^{\frac 12}} dt
\leq  \log (3+\eps)-\log 3.
\]
Thus  a straightforward computation gives
\[
\int_0^R\frac{\sigma\sigma'\sigma''}
{(1+\sigma'(t)^2)^{\frac 32}}
dt\leq
-\frac{\log R}{(1+\frac{1}{R^2})^\frac 12}
+\log( 3+\eps) +\log (R+\sqrt{1+R^2})-\log (3+\eps
+\sqrt{1+(3+\eps)^2}).
\]
We let $R\rightarrow \infty, \eps\rightarrow 0$.
Then we have
\begin{equation}\label{ex-n-5}
\int_0^\infty\frac{\sigma\sigma'\sigma''}
{(1+\sigma'(t)^2)^{\frac 32}}
dt\leq\log 6-\log (3+\sqrt{10})<0.
\end{equation}
By~\eqref{ex-n-3} and~\eqref{ex-n-5}, we have
$\int_M c_2(A)<0$. 
Finally, since 
\[\frac{\sigma\sigma'\sigma''}
{(1+\sigma'(t)^2)^{\frac 32}}\sim
O(\frac{\log t}{t^3}),
\]
and 
\[
-\frac{t\sigma''}{(1+\sigma'(t)^2)^\frac 32}
-\frac{\sigma'}{\sqrt{1+\sigma'(t)^2}}
\sim O(\frac{1}{t^3}),
\]
we know that $c_2(A)$ is integrable.

\qed

\section{Convex surfaces}

In this section, we consider the layer over a convex surface
$\Sigma$ in
$\R^3$. $\Sigma$ is defined by 
\begin{equation}
z=f(x,y),
\end{equation}
where $f(x,y)$ is a smooth convex function,  $f(0)=0$,
$\nabla f(0)=0$, and $\nabla^2 f(0)>0$.

The main result of this section is the following:

\begin{theorem}\label{thm61} Let $\Sigma$ be defined as above.
 Suppose
$\Omega=\Sigma\times(-a,a)$ is the layer with depth $a>0$.
Then the infimum of the spectrum $\sigma_0$ satisfies
\begin{equation}
\sigma_0<\kappa_1^2,
\end{equation}
where $\kappa_1=\pi/(2a)$.
\end{theorem}

We begin with the following 

\begin{lemma} 
\label{lem61}
With the  assumptions on $f$,  there
is a number $\delta>0$ such that
\[
f_{r}=\frac{\pa f}{\pa r}>\delta
\]
for $x^2+y^2\geq 1$, where $(r,\theta)$ is the polar
coordinates defined by $x=r\cos\theta$ and $y=r\sin\theta$.
\end{lemma}

{\bf Proof.} 
By the assumptions, there is a $\delta_0\in(0,1)$ such that
$f_{r}|_{x^2+y^2=\delta_0}>0$. By convexity, we have
$f_{rr}\geq 0$. Then since the circle $\{x^2+y^2=\delta_0\}$ is 
compact, we can conclude that $f_r\geq\delta$ for $x^2+y^2\geq 1$.

\qed

\begin{cor} \label{cor61}
Using the above notations, we have 
\begin{enumerate}
\item $|\nabla f|\geq\delta$;
\item $f(x,y)\geq \delta\cdot (\sqrt{x^2+y^2}-1)$.
\end{enumerate}
\end{cor}

\qed

An interesting consequence of the above corollary is the
following: Let $b$ be a large positive number. Let $C_b$ be
the curve defined by the intersection of $\Sigma$ with
respect to the  plane $z=b$. Clearly $C_b$ is a convex
curve. From the above corollary, $C_b$ is contained 
in a disk of radius $b/\delta+1$. In particular, we
have the estimate of the length of the curve
\begin{equation}\label{6-2}
\int_{C_b}1\leq Cb
\end{equation}
for a constant $C$.

For the manifold $\Sigma$, the mean curvature $H$ can be
represented by
\begin{equation}
H=\frac{(1+f_y^2)f_{xx}+(1+f_x^2)f_{yy}-2f_xf_yf_{xy}}
{(1+|\nabla f|^2)^{3/2}},
\end{equation}
where  $f_x=\frac{\pa f}{\pa x}$, 
$f_y=\frac{\pa f}{\pa y}$, $f_{xx}=
\frac{\pa^2 f}{\pa x^2}$, $f_{yy}=
\frac{\pa^2 f}{\pa y^2}$, $f_{xy}=
\frac{\pa^2 f}{\pa x\pa y}$, and $|\nabla f|^2=f_x^2
+f_y^2$.
We  compare the mean curvature  to 
the curvature of the convex curve $f(x,y)=b$, which 
is given by
\begin{equation}
k_b=\frac{f_{xx}f_y^2-2f_{xy}f_xf_y+f_{yy}f_x^2}{|\nabla
f|^{3}}.
\end{equation}
By Corollary~\ref{cor61} and the convexity of $f$, we have
\begin{equation}\label{6-4}
H\geq \frac 12\delta^3 k_b,
\end{equation}
if $\delta$ is small enough.
Since $C_b:\{f(x,y)=b\}$ is a convex curve, we have
\begin{equation}
\int_{f(x,y)=b}k_b=2\pi.
\end{equation}
Thus by~\eqref{6-4}
\begin{equation}\label{6-7}
\int_{f(x,y)=b}H\geq\pi\delta^3. 
\end{equation}

By the co-area formula (cf.~\cite[page 89]{SY}), we have

\begin{equation}
\int_{x^2+y^2\geq 1} Hd\Sigma\geq\int^\infty_{c}
(\int_{f=t}\frac{H}{|\tilde \nabla f|})
dt,
\end{equation}
where $c$ is a positive real  number, and
\begin{equation}\label{6-9}
|\tilde\nabla f|^2=\frac{|\nabla f|^2}
{1+|\nabla f|^2}
\end{equation}
is the gradient of $f$ on the {\it Riemannian manifold}
$\Sigma$. Thus by ~\eqref{6-7},
and Corollary~\ref{cor61},
\begin{equation}
\int_{x^2+y^2\geq 1} H=+\infty.
\end{equation}

{\bf Proof of Theorem~\ref{thm61}.} We shall 
again use the  trick
introduced by~\cite{DEK} (see also ~\cite{DP,RB}) to  
perturb the 
``standard'' test functions. However,
our choices of perturbation functions are quite  
different from  theirs in nature.

Let $K$ be the Gauss curvature of $\Sigma$. Then we have
$K\geq 0$, and 
\begin{equation}\label{6-11}
\int_\Sigma K\leq 2\pi
\end{equation}
by the theorem of Huber~\cite{Huber}.
Since the Gauss curvature is nonnegative, the volume
growth is at most quadratic. Thus
$\Sigma$ is parabolic. For any $r_1>0$, we can find a
function
$\phi$ such that 
\begin{enumerate}
\item $\phi\in C_0^\infty(\Sigma)$,
$0\leq\phi\leq 1$;
\item $\phi\equiv 1$ on $B(r_1)$,
where $B(r_1)$ is the geodesic ball of radius $r_1$
of $\Sigma$
centered at $0$;
\item $\int_\Sigma|\nabla\phi|^2d\Sigma<1$.
\end{enumerate}

The quadratic forms $Q$, $Q_1$, and $Q_2$ are 
defined in~\eqref{5-1-1},~\eqref{5-2-1}, and~\eqref{5-3}.
Let $\chi=\cos\kappa_1 u$. Then we  have
\begin{equation}
Q_1(\phi\chi,\phi\chi)=\int_\Omega|\nabla\phi|^2\chi^2d\Omega
\leq a{(1+C_0)^2},
\end{equation}
where $C_0<1$ is defined in Assumption A1).
We also have
\begin{equation}
Q_2(\phi\chi,\phi\chi)\leq
\mu_2\int_{\Sigma}K\phi^2d\Sigma.
\end{equation}
Combining the above two equations and using~\eqref{6-11},
we  have
\begin{equation}\label{6-12}
Q(\phi\chi,\phi\chi)=Q_1(\phi\chi,\phi\chi)
+Q_2(\phi\chi,\phi\chi)\leq C_1,
\end{equation}
where $C_1$ is a constant depending only on $\Sigma$
and $a$.

Suppose $r_1$ is large enough such that
$\{f(x,y)\leq 2R^2\}\subset B(r_1)$ for some 
large number $R>0$.
 We consider a function
$\rho(t)$ on $\R$ such that 
\begin{enumerate}
\item $\rho\equiv 1$, if $t\in[R,R^2]$;
\item $\rho \equiv 0$ if $t>R^2+1$ or $t<R-1$;
\item $0\leq\rho\leq 1$;
\item $|\rho'|\leq 4$.
\end{enumerate}

We define $\psi(x,y)=\rho(f(x,y))$. $\psi$ is a smooth
function of $\Sigma$.
Let $\chi_1$ be an odd function of $u$
such that $\chi_1(\pm a)=0$, and
\begin{equation}\label{6-15}
\int_{-a}^a\chi_u(\chi_1)du=-\sigma<0,
\end{equation}
where $\sigma>0$ is a positive number.
 We consider
the function
$\phi\chi+\eps\psi\chi_1/f$, where $\eps$ is a small number
to be determined. By the definition of $Q(\cdot,\cdot)$
and~\eqref{6-12}, we
have
\begin{equation}\label{6-16}
Q(\phi\chi+\eps\psi\chi_1/f, \phi\chi+\eps\psi\chi_1/f)
\leq
C_1+2\eps Q(\phi\chi,\psi\chi_1/f)+\eps^2
Q(\psi\chi_1/f,\psi\chi_1/f).
\end{equation}
If $r_1$ and $R$ are big, then 
\[
{\rm supp}\,\psi\subset \{x\in\Sigma\mid \phi(x)\equiv 1\}.
\]
We thus have
\[
\int_\Omega\langle \nabla' (\phi\chi),\nabla'(\psi\chi_1/f)\rangle d\Omega=0,
\]
where $\nabla'$ is defined in~\eqref{5-2-2}.
By~\eqref{2-3}, we have
\[
d\Omega=(1-Hu+Ku^2)d\Sigma.
\]
Since $\chi_1$ is odd and $\chi$ is even, by the above equation, we have
\[
Q(\phi\chi,\psi\chi_1/f)=\int_\Omega(\chi_u(\chi_1)_u\phi\psi/f-\kappa_1^2
\chi\chi_1\phi\psi/f)d\Omega
=-\int^a_{-a}u(\chi_u(\chi_1)_u-\kappa_1^2\chi\chi_1) du\int_\Sigma
\psi H/fd\Sigma.
\]
Since 
\[
\int^a_{-a}u(\chi_u(\chi_1)_u-\kappa_1^2\chi\chi_1) du=-\int_{-a}^a\chi_u(\chi_1)du,
\]
 we have
\begin{equation}\label{6-17}
Q(\phi\chi,\psi\chi_1/f)=-\sigma\int_\Sigma
\frac{\psi H}{f}d\Sigma,
\end{equation}
where $\sigma$ is the number defined in ~\eqref{6-15}.
By the co-area formula and ~\eqref{6-7}, ~\eqref{6-9}, we
have
\begin{equation}\label{6-18}
\int_\Sigma\frac{\psi H}{f}=\int_{\mathbb R}\frac
{\rho(t)}{t}
(\int_{f=t}\frac{H}{|\tilde\nabla
f|})dt\geq\pi\delta^3\int_{\mathbb R}\frac{\rho(t)} t
dt\geq\pi\delta^3\log R.
\end{equation}

In order to estimate the last term of ~\eqref{6-16}, we
first note that
\begin{equation}
|\tilde \nabla\frac\psi f|\leq C_2/f
\end{equation}
for some constant $C_2$, where $\tilde\nabla$ is the 
covariant derivative of $\Sigma$. From this,
 we have
\begin{equation}\label{6-20}
Q_1(\frac{\psi\chi_1}{f},\frac{\psi\chi_1}{f})\leq
C_3\int_{R-1\leq f\leq R^2+1} 1/f^2d\Sigma
\end{equation}
for some constant $C_3$.
Using the same argument, we have
\begin{equation}\label{6-21}
Q_2(\frac{\psi\chi}{f},\frac{\psi\chi}{f})
\leq
C_4\int_{R-1\leq f\leq R^2+1} 1/f^2d\Sigma
\end{equation}
for some constant $C_4$.
We use the co-area formula again to estimate
\begin{equation}
\int_{R-1\leq f\leq R^2+1} 1/f^2
=\int_{R-1}^{R^2+1}\frac{1}{t^2}(\int_{f=t}\frac{1}{|\tilde
\nabla f|})dt.
\end{equation}
From Corollary~\ref{cor61}, we know that
$|\tilde\nabla f|$ has a lower bound. Thus 
by~\eqref{6-2}, there is a constant
$C_5$ such that 
\begin{equation}\label{6-23}
\int_{R-1\leq f\leq R^2+1} 1/f^2d\Sigma\leq
C_5\log R.
\end{equation}
Thus by~\eqref{6-17},~\eqref{6-18},~\eqref{6-20}, 
~\eqref{6-21}, and~\eqref{6-23}, from ~\eqref{6-16}, we
have
\begin{equation}\label{6-25}
Q(\phi\chi+\eps\psi\chi_1/f, \phi\chi+\eps\psi\chi_1/f)
\leq C_1 -2\eps\pi\sigma\delta^3\log
R+\eps^2C_5(C_3+C_4)\log R.
\end{equation}
We choose $\eps$ to be a small positive number
such that
\[
-2\eps\pi\sigma\delta^3+\eps^2C_5(C_3+C_4)<0.
\]
We then let $R$ large enough (which requires $r_1$ be large
enough also). Then the left hand side of ~\eqref{6-25} is
negative. By the definition of $\sigma_0$, we know that
$\sigma_0<\kappa_1^2$. 

\qed

{\bf Proof of Theorem~\ref{main3}.}
By Theorem~\ref{thm31} and Theorem~\ref{thm61}, we have
\[
\sigma_0<\kappa_1^2\leq\sigma_{\rm ess}.
\]
Thus the ground state exists.

\qed

\bibliographystyle{abbrv}   
\bibliography{new060923,unp060923,local}

\end{document}